\documentclass[12pt,a4paper]{amsart}

\usepackage{amssymb}  

\newcommand{\Z}{{\mathbb Z}}
\newcommand{\Q}{{\mathbb Q}}

\newcommand{\BP}{{\mathbb P}}
\newcommand{\CO}{{\mathcal O}}
\newcommand{\CC}{{\mathcal C}}
\newcommand{\CD}{{\mathcal D}}
\newcommand{\CK}{{\mathcal K}}
\newcommand{\CW}{{\mathcal W}}
\newcommand{\CX}{{\mathcal X}}
\newcommand{\Fp}{{\mathfrak p}}
\newcommand{\To}{\longrightarrow}
\newcommand{\tensor}{\otimes}
\newcommand{\disc}{\operatorname{disc}}
\newcommand{\Aut}{\operatorname{Aut}}
\newcommand{\Spec}{\operatorname{Spec}}
\newcommand{\Pic}{\operatorname{Pic}}
\newcommand{\Sym}{\operatorname{Sym}}
\newcommand{\half}{\tfrac{1}{2}}
\newcommand{\inj}{\hookrightarrow}
\newcommand{\unr}{{\text{\rm unr}}}
\newcommand{\s}{{\text{\rm s}}}
\newcommand{\tors}{{\text{\rm tors}}}
\newcommand{\triv}{{\text{\rm $\Gamma$-triv}}}
\newcommand{\Triv}{{\text{\rm triv}}}
\newcommand{\ceil}[1]{\left\lceil #1 \right\rceil}
\newcommand{\codim}{\operatorname{codim}}

\newtheorem{Theorem}{Theorem}[section]
\newtheorem{Lemma}[Theorem]{Lemma}
\newtheorem{Proposition}[Theorem]{Proposition}
\newtheorem{Corollary}[Theorem]{Corollary}

\theoremstyle{remark}
\newtheorem{Remark}[Theorem]{Remark}

\numberwithin{equation}{section}

\begin{document}

\title[Independence of rational points]%
      {Independence of rational points\\
       on twists of a given curve}

\author{Michael Stoll}
\address{School of Engineering and Science,
         International University Bremen,
         P.O.Box 750561,
	 28725 Bremen, Germany.}
\email{m.stoll@iu-bremen.de}
\date{February 27, 2006}


\begin{abstract}
  \setlength{\parskip}{0.4ex plus 0.1ex minus 0.1ex}
  \setlength{\parindent}{0mm}
  
  In this paper, we study bounds for the number of rational points
  on twists~$C'$ of a fixed curve~$C$ over a number field~$\CK$, under
  the condition that the group of $\CK$-rational points on the
  Jacobian~$J'$ of~$C'$ has rank smaller than the genus of~$C'$.
  
  The main result is that with some explicitly given finitely many possible 
  exceptions, we have a bound of the form $2r + c$, where $r$ is the
  rank of~$J'(\CK)$ and $c$ is a constant depending on~$C$.
  
  For the proof, we use a refinement of the method of Chabauty-Coleman;
  the main new ingredient is to use it for an extension field of~$\CK_v$,
  where $v$ is a place of bad reduction for~$C'$.
\end{abstract}

\maketitle


\section{Introduction}

Let $C/\CK$ be a smooth projective curve over a number field. 
Fix some $\CK$-rational
divisor class $\CD$ of positive degree~$d$ (for example, the canonical
class, if the genus of~$C$ is at least~$2$) and use it as a ``basepoint''
in order to map points on~$C$ to points on the Jacobian~$J$ of~$C$:
\[ \phi : C \ni P \longmapsto [d \cdot P] - \CD \in J \]
Now, given a set $\Sigma \subset C(\CK)$ of rational points on~$C$, we can
ask ourselves how large the subgroup of~$J(\CK)$ generated by~$\phi(\Sigma)$
may be. In particular, are the points in~$\Sigma$ indpendent, i.e., does
$\phi(\Sigma)$ generate a group of rank~$\#\Sigma$?

We can also turn around this question --- we take a subgroup $G \subset J(\CK)$
and ask how many points in~$C(\CK)$ map into~$G$. In this paper, we will
give answers to these questions when $C$ is a twist of a fixed curve.
It will turn out that we get fairly tight bounds if the number of points
(or the rank of the subgroup) is sufficiently small relative to the genus,
as long as we are willing to accept finitely many exceptions (which can
be found in an explicitly given finite set of twists).

We give a few applications. The first one deals with quadratic twists
of hyperelliptic curves. For simplicity, we formulate the result with
$\Q$ as the base field, though it is valid for any number field.
For hyperelliptic curves, we use the class of twice a Weierstrass point
as our ``basepoint''~$\CD$.

\begin{Theorem} \label{IntroThmHC}
  Let $C : y^2 = f(x)$ be a hyperelliptic curve over~$\Q$ of genus $g \ge 2$, 
  where $f \in \Z[x]$ is squarefree. Let $\Delta$ be the discriminant
  of~$f$, considered as a polynomial of degree~$2g+2$.
  Let $d \in \Z$ be a squarefree integer
  and $n \le g$ a natural number such that $d$ is divisible by a 
  prime~$p > 2n + 1$ that does not divide~$\Delta$.
  Let $C_d : d\,y^2 = f(x)$ be the
  quadratic twist of~$C$ associated to~$d$, and denote by~$\iota$ the
  hyperelliptic involution. Then any set $\Sigma \subset C_d(\Q)$ of rational
  points on~$C_d$ such that $\#\Sigma \le n$ and 
  $\Sigma \cap \iota(\Sigma) = \emptyset$
  generates a subgroup of rank~$\#\Sigma$ in the Jacobian of~$C_d$.
\end{Theorem}

Note that the exceptional values of~$d$ are squarefree integers with prime 
divisors
in the finite set \hbox{$\{p \mid p \le 2n+3 \text{\ or\ } p | \Delta\}$}
and so are finite in number. The condition 
$\Sigma \cap \iota(\Sigma) = \emptyset$ is necessary, since 
$\phi(\iota(P)) = -\phi(P)$ provides a trivial dependence.

Seen from the other side, we can state this result in the following form.

\begin{Theorem} \label{IntroThmHCVariant}
  Keep the notations of Thm.~\ref{IntroThmHC}. Assume that the Mordell-Weil
  rank $r$ of the Jacobian of~$C_d$ satisfies $r < g$ and that the prime~$p$
  dividing~$d$ satisfies $p > 2r + 3$ and does not divide~$\Delta$.
  Then $C_d$ has at most $r$ pairs of rational non-Weierstrass points.
\end{Theorem}

Again there are only finitely many exceptional squarefree~$d$, 
contained in an explicit finite set. This result is obviously best possible,
except that it does not say anything about the excluded cases.

The next application is to Thue equations.

\begin{Theorem} \label{IntroThmThue}
  Let $F \in \Z[X,Y]$ be homogeneous of degree~$n \ge 3$ and squarefree.
  Let $h \in \Z$ be an integer not divisible by the $n$th power of any prime,
  such that $h$ has a prime factor $p > n+r+1$ that does not divide the
  discriminant of~$F$, where $r$ is the Mordell-Weil rank of the Jacobian
  of the curve given by the Thue equation
  \begin{equation} \label{EqnThue}
     F(X, Y) = h \,. 
  \end{equation}
  If $r \le n-3$, then this equation has at most~$r$ rational solutions.
  More generally, if $r \le \half n(n-3)$ and $p > n + 2r + 1$, 
  then there are at most~$2r$ rational solutions.
\end{Theorem}

The possible exceptions are again finite in number and contained in an
explicitly given set. It is interesting to compare this with the bound
of Lorenzini and Tucker \cite{LorenziniTucker}, Thm.~3.11, which is
that there are at most $2n^3 - 2n - 3$ primitive {\em integral} solutions
to~\eqref{EqnThue} if $r \le \half n(n-3)$. The bound is weaker, but
it holds for all~$h$ (subject to the rank condition). On the other hand,
our result, when applicable, even bounds the number of {\em rational}
solutions, and our bound is {\em much} stronger.

In a similar way, we can state a general bound on the number of rational points
on twists of a fixed curve. For a curve~$C$, denote by $C_\Triv$ the
set of points that are fixed by some nontrivial (geometric) automorphism
of~$C$ and by~$C_\tors$ the set of points that map to a torsion point
in the Jacobian, where we have chosen the ``basepoint''~$\CD$ to be
invariant under the automorophism group of~$C$. We denote by $r(C)$
the Mordell-Weil rank of $J(\CK)$, where $J$ is the Jacobian of~$C$.

\begin{Theorem} \label{IntroThmPointsOnTwists}
  Fix a curve $C$ of genus~$g \ge 2$ over a number field~$\CK$. 
  Then for all but finitely many twists $C'/\CK$ of~$C$ such that $r(C') < g$,
  we have
  \[ \#C'(\CK) \le f_C(r(C')) + \#C'_\Triv(\CK)
                    + \#(C'_\Triv \setminus (C'_\Triv(\CK) \cup C'_\tors)) \,.
  \]
  Here, $f_C$ is a function that depends on the geometry of~$C$; for 
  $0 \le r < g$, we have $r \le f_C(r) \le 2r$.
\end{Theorem}

In particular, this implies that within the family of twists of rank $< g$
of~$C$, the number of rational points is bounded. This is not new;
it follows from a result due to Silverman~\cite{Silverman93},
which states that 
\[ \#C'(\CK) \le \gamma(C/\CK) \, 7^{r(C')} \]
for all twists~$C'$, with a constant $\gamma(C/\CK)$ depending only on~$C/\CK$,
which is effective in principle, but which nobody has tried to find a value
for, as far as we know. Silverman's result is stronger than ours in that
it covers all twists, regardless of the rank. On the other hand, our result
is stronger than Silverman's in that it provides a much better bound, when
it applies. It should be mentioned, however, that the so-called Lang conjecture
on varieties of general type would imply that the number of rational points
on any curve (of genus at least~$2$) is bounded by a constant only depending
on the genus and the degree of the base field, see
\cite{CHM} and~\cite{Pacelli}.

As a final application, we remark that one can use the method developed
here to obtain sharp bounds for the number of rational points on twists
of fixed rank. For example, a detailed study of the possible behaviour
at small primes shows that for the family
\[ C_A : y^2 = x^5 + A \]
(with $A$ an integer not divisible by the tenth power of any prime), the
maximal number of rational points is~$7$ when $r(C_A) = 1$, and this
maximum is attained only for the curve that has $A = 18^2$. Note that
the general result proved in this paper implies that all but finitely
many of these curves (such that $r(C_A) = 1$) have at most $5$ rational
points. The more specific arguments necessary to obtain the precise
result stated above are given in~\cite{Stoll06}.

As is already apparent to the cognoscenti from the small rank conditions
in the results given above, we use
a version of the Chabauty-Coleman method for the proof. 
The main new ingredient
is to apply this method over an {\em extension field} of a completion~$\CK_v$,
where $v$ is a place of {\em bad reduction} for the twist in question. 
Surprisingly, this
leads to much better bounds than working directly with~$\CK_v$ or at
primes of good reduction.

In the next sections, we state the main theorem and deduce the results
given above. We then proceed to review the Chabauty-Coleman method.
Finally, we apply this method to prove our main theorem.

\subsection*{Acknowledgements}

First of all, I wish to thank MSRI for inviting me to spend a month
there (mid-November to mid-December 2000). The results described here
have their orgin in work I did during this stay. 
Then, but certainly not less important,
I want to thank Nils Bruin, Noam Elkies, Dino Lorenzini, Bjorn Poonen,
Ed Schaefer, Joe Silverman and Joe Wetherell for useful, interesting 
and sometimes inspiring conversations, either directly or by email.


\section{Notation}

For a field~$\CK$, we denote by $\bar{\CK}$ a fixed algebraic closure.
When $\Gamma$ is a (not necessarily abelian) group, on which the
absolute Galois group of~$\CK$ acts, we denote by $H^1(\CK, \Gamma)$
the first Galois cohomology; for non-abelian $\Gamma$, this is a 
pointed set (see~\cite{Serre}, \S~5).
We let the Galois group act on the right and write the action 
exponentially: $P^\sigma$.

We continue to use the notation from the introduction. So $C/\CK$
is a smooth projective curve over a number field.
The Jacobian of~$C$ is denoted $J$. It is an abelian variety over~$\CK$
of dimension~$g = g(C)$, the genus of~$C$. We fix a $\CK$-rational divisor
(or divisor class) $\CD$ on~$C$ of degree $d > 0$ (which will have
to satisfy a certain invariance condition, see below; if $g \ge 2$,
we  can always take the canonical class) and use it to map
$C$ into~$J$ via $\phi : P \mapsto [dP - \CD]$. 

It then makes sense
to define $C_\tors$ as the preimage of the torsion points in~$J$.
Note that this set is finite when the genus of~$C$ is at least two,
see Raynaud's proof of the Manin-Mumford conjecture in~\cite{Raynaud}.


\section{Some Geometry} \label{SectGeom}

For a $\bar{\CK}$-defined divisor $D$ on~$C$, we let $\Omega(D)$ denote the 
$\bar{\CK}$-vector space
of differentials~$\omega$ satisfying $(\omega) \ge D$. We then define
the function $f_C : \Z_{\ge 0} \to \Z_{\ge 0} \cup \{\infty\}$ as follows.
\[ f_C(r)
    = \max\{\deg D \mid D \ge 0 \,\text{\ and\ } \dim \Omega(D) \ge g - r\} \,,
\]
where $D$ runs through all effective divisors on~$C$.
We obviously have $f_C(r) = \infty$ as soon as $r \ge g$.
For the interesting values of~$r$, we have the following results.

\begin{Lemma} \label{LemmaFBound} \strut
  \begin{enumerate}
    \item \label{S1} If $0 \le r < g$, then $r \le f_C(r) \le 2r$. 
    \item \label{S2} $f_C(0) = 0$, $f_C(g-1) = 2g-2$.
    \item \label{S3} $C$ is hyperelliptic if and only if $f_C(1) = 2$, if and
          only if $f_C(r) = 2r$ for $0 \le r < g$.
    \item \label{S4} If $C$ is a smooth plane curve of degree~$n$, then
      \[ f_C(r) = r + \binom{n-a}{2} - 1 \,\text{, where\ }
              a = \max \Bigl\{k \Bigm|
                              r + \binom{k}{2} < g = \binom{n-1}{2} \Bigr\} \,. 
      \]
          In particular, $f_C(r) = r$ for $r \le n-3$ in this case.
  \end{enumerate}
\end{Lemma}
\begin{proof}
  \ref{S1}) The lower bound follows from $\dim \Omega(D) \ge g - \deg D$ for any
      effective divisor~$D$. For the upper bound, we use 
      Riemann-Roch and the standard bound
      $\dim L(D) \le 1 + \deg D/2$ for $D$ such that $0 \le \deg D \le 2g$.
      This gives
      \[ g - r \le \dim \Omega(D) = \dim L(D) - \deg D + g - 1 \le g - \deg D/2        \]
      and so $\deg D \le 2r$.
  \ref{S2}) Clear.
  \ref{S3}), \ref{S4}) The facts about hyperelliptic and plane curves are well-known.
\end{proof}

When we consider $f_{C/k_v}(r)$ for the reduction of $C$ at a place
$v$ of~$\CK$ of good reduction for~$C$ (where $k_v$ is the residue
field), it is possible that $f_{C/k_v}(r) > f_C(r)$ for certain $v$ and~$r$.
For example, the Klein Quartic, a smooth plane quartic curve, 
has over a certain number field~$\CK$ that is
totally ramified at~$7$ a model that reduces to a hyperelliptic curve
at the place above~$7$, see~\cite{Elkies}. In this case, we have
$f_C(1) = 1$, but $f_{C/k_v}(1) = 2$.

However, we have the following result.

\begin{Proposition} \label{PropBadFinite}
  For all but finitely many places $v$ of~$\CK$ such that $C$  has
  good reduction at~$v$, we have $f_{C/k_v} = f_C$.
\end{Proposition}
\begin{proof}
  Let
  \[ F_C(d) = \max\{\dim L(D) \mid D \ge 0 \text{\ and\ } \deg D = d\} \]
  be the largest dimension of a Riemann-Roch space of an effective
  divisor of degree~$d$; define $F_{C/k_v}$ similarly. Then, since
  by the Riemann-Roch Theorem,
  \[ \dim L(D) = \deg D - g + 1 + \dim \Omega(D) \,, \]
  we have, for $0 \le r < g$,
  \[ f_C(r) = \max\{d \mid F_C(d) \ge d - r + 1\} \,, \]
  and similary for $f_{C/k_v}$. Therefore it is sufficient to show
  the statement with $F_C$ and $F_{C/k_v}$ in place of $f_C$
  and~$f_{C/k_v}$.
  
  Now choose a smooth projective model $\CC/\CO_S$ of $C$ over
  some ring of $S$-integers $\CO_S$ in~$\CK$. Then we have a
  canonical morphism
  \[ \phi : \Sym_{\CO_S}^d \CC \To \Pic_{\CO_S}^d \CC =: \CW \]
  from the $d$th symmetric power of the $\CO_S$-scheme $\CC$
  into the degree-$d$ component of the relative Picard scheme
  of $\CC$ over~$\CO_S$. According to our definition,
  \[ F_C(d) = 1 + \max\{\dim \phi^{-1}(w) \mid w \in \CW(\bar{\CK})\} \]
  and
  \[ F_{C/k_v}(d) = 1 + \max\{\dim \phi^{-1}(w) 
                                \mid w \in \CW(\bar{k}_v)\} \,.
  \]
  Now the subscheme of~$\CW$ of points~$w$ such that the dimension
  of the fiber $\phi^{-1}(w)$ is at least a certain number~$n$ is constructible
  and so is its image $\CX_n$ in~$\Spec \CO_S$ under the structure morphism
  (see \cite{Hartshorne}, Exc.~II.3.19 and~II.3.22, or
  \cite{CartanChevalley}, exp.~7 and~8). Also, these schemes will
  be empty for $n > d-g$. For $0 \le n \le F_C(d)$, $\CX_n$ will
  contain the generic point of~$\Spec \CO_S$; for larger $n$, it will not.
  The set of places~$v$ outside~$S$ such that $F_{C/k_v}(d) \neq F_C(d)$
  is therefore contained in the constructible subscheme
  \[ \CX(d) = \bigcup_{n=0}^{F_C(d)} (\Spec \CO_S \setminus \CX_n)
                 \cup \bigcup_{n = F_C(d)+1}^{d-g} \CX_n \,.
  \]
  Since $\CX(d)$ does not contain the generic point, it is finite,
  and since we have to consider only $d \le 2g-2$, the set of
  places~$v$ such that $F_{C/k_v} \neq F_C$ is also finite.
\end{proof}

It is a different matter to actually {\em find} the ``mildly bad''
places~$v$ such that the curve has good reduction at~$v$,
but ``its geometry changes'' under reduction, in the sense that
$f_C$ changes. For practical purposes, we will therefore choose
a function $\tilde{f}_C$, together with a set of ``bad'' places
$S_C$ such that for all places $v \notin S_C$, $C$ has good reduction
at~$v$ and $f_{C/k_v} \le \tilde{f}_C$. Of course, Prop.~\ref{PropBadFinite}
says that there is always some such set~$S_C$ that works with
$\tilde{f}_C = f_C$.

\begin{Proposition} \label{PropChooseF}
  The following are possible choices for $\tilde{f}_C$ and~$S_C$.
  \begin{enumerate}
    \item $\tilde{f}_C(r) = 2r$ and $S_C =$ places of bad reduction.
    \item If $C$ is a smooth plane curve, pick a projective plane
          model $\CC/\CO$ over the integers of~$\CK$ and take 
          $\tilde{f}_C = f_C$ and $S_C =$ places of bad reduction
          of the model~$\CC$.
  \end{enumerate}
\end{Proposition}
\begin{proof}
  Clear.
\end{proof}

Note that hyperelliptic curves are always ``worst-case'' for $f_C$,
and so we lose nothing if we take the first alternative above in this case.

\bigskip

For the following, we let $\CK$ be an arbitrary field.
Recall that the set of {\em twists} of~$C$, i.e., the set of curves
$C'/\CK$, isomorphic to~$C$ over~$\bar{\CK}$, up to isomorphism
over~$\CK$, is parametrized by the Galois cohomology set
$H^1(\CK, \Aut_{\bar{\CK}}(C))$. Given such a twist $C'$ and a
$\bar{\CK}$-isomorphism $\varphi : C' \To C$, the corresponding
cohomology class is represented by the cocycle
$\xi : \sigma \mapsto \varphi^\sigma \circ \varphi^{-1}$; see for
example Silverman's book~\cite{SilvermanBook}, \S~X.2.

For the purposes of this paper, it is useful to introduce a more
general notion. A twist $C'$ of~$C$ is called a {\em $\Gamma$-twist} 
of~$C$ if the corresponding cohomology class lies in the image
of~$H^1(\CK, \Gamma)$ in~$H^1(\CK, \Aut_{\bar{\CK}}(C))$. This means
that there is a \hbox{$\bar{\CK}$-isomorphism} $\varphi : C' \To C$ such
that the cocycle $\xi : \sigma \mapsto \varphi^\sigma \circ \varphi^{-1}$
takes values in~$\Gamma$. For such a cocycle, we will denote the
corresponding $\Gamma$-twist of~$C$ by~$C_\xi$.

An important example is given by quadratic twists of hyperelliptic
curves; here we take for $\Gamma$ the subgroup of order two of the
automorphism group generated by the hyperelliptic involution.


\section{The main result --- local version}

In this section, we consider a smooth projective
curve $C$ of genus $g \ge 1$ over a $p$-adic field~$K$.
We denote by~$v$ its normalized additive valuation, by~$\CO$
its ring of integers, and by~$k$ its residue field. $G_K$ is the
absolute Galois group of~$K$, relative to a fixed algebraic 
closure~$\bar{K}$, and $I_K \subset G_K$ denotes the inertia
group. $K^{\unr}$ is the maximal unramified extension
of~$K$ (inside~$\bar{K}$).

We will always assume that $C$ has good reduction, so that there
is a smooth projective curve $\CC/\CO$ with generic fiber~$C$.

Let $\Gamma \subset \Aut_{\bar{K}}(C)$ be a finite $K$-defined
subgroup of the (geometric) automorphism group of~$C$. We always
assume that $v(\#\Gamma) = 0$, i.e., that the order of~$\Gamma$
is prime to the residue characteristic~$p$.

\begin{Proposition} \label{PropIActsTriv}
  Under the assumptions made, $I_K$ acts trivially on~$\Gamma$.
\end{Proposition}

\begin{proof}
  Let $\gamma \in \Gamma$ and $\sigma \in I_K$, and set
  $\phi = \gamma^{-1} \gamma^\sigma \in \Gamma$. We have to
  show that $\phi$ is the identity. Now, since $g \ge 1$ and
  $\CC$ is the minimal proper regular model of~$C$ (trivially),
  $\phi$ extends to an automorphism of~$\CC$, see for example
  \cite{SilvermanBook2}, Prop.~IV.4.6. (We base-extend $\CC$
  to the ring of integers of the field of definition of~$\phi$;
  since we have good reduction, we still have a minimal proper
  regular model.) Since $I_K$ acts trivially on the special fiber
  of~$\CC$, the automorphism induced by~$\phi$ on the special
  fiber is the identity. Pick some point $P \in \CC(\bar{k})$
  and consider its residue class
  \[ D_P = \{Q \in C(\bar{K}) \mid \bar{Q} = P\} \,, \]
  (where $\bar{Q}$ is the image of $Q$ in the special fiber of~$\CC$).
  We claim that either $\phi$ is the identity on~$D_P$ or else
  it has a unique fixed point on~$D_P$. Given that, it follows 
  that $\phi = 1$ in~$\Gamma$: in any case, $\phi$ has
  at least one fixed point on every residue class $D_P$; since
  there are infinitely many residue classes, $\phi$ has infinitely
  many fixed points on~$C$ and thus must be the identity automorphism.
  
  Let us now prove the claim. Consider a finite extension
  $L$ of~$K$ such that $\phi$ is defined over~$L$ and $P$ is defined
  over the residue class field of~$L$. Then
  \[ D_P(L) = D_P \cap C(L) \]
  can be parametrized analytically by the maximal ideal $\Fp_L$
  of~$\CO_L$, and in this parametrization, the action of~$\phi$
  is given by a power series
  \[ F_\phi(T) = \alpha_0 + \alpha_1 T + \dots \in \CO_L[T] \]
  with $v_L(\alpha_0) > 0$, since the action fixes~$\Fp_L$.
  Now $\phi^m = 1$ with some $m$ prime to~$p$, so, reducing mod~$\Fp_L$,
  \[ \bar{F}_\phi^{\circ m}(T) = \bar{\alpha}_1^m T + \dots = T \,. \]
  Hence $\alpha_1$ is a unit reducing to an $m$th root of unity. 
  Suppose first that $\bar{\alpha}_1 = 1$ and that $F_\phi(T) \neq T$
  (otherwise, $\phi$ is the identity on~$D_P$). Write
  \[ F_\phi(T) = T + \pi_L^n \tilde{F}(T) \]
  for the maximal possible~$n$; we claim that $n \ge 1$. Otherwise,
  write (mod~$\Fp_L$)
  \[ \bar{F}_\phi(T) = T + \beta T^N + O(T^{N+1}) \]
  with $N \ge 2$ and $\beta \neq 0$; then
  \[ T = \bar{F}_\phi^{\circ m}(T) = T + m \beta T^N + O(T^{N+1}) \,, \]
  a contradiction. So in the above, $n \ge 1$, and
  $\tilde{F}$ will not reduce to zero mod~$\Fp_L$. Iterating, we get
  \[ T = F_\phi^{\circ m}(T) = T + m\,\pi_L^n \tilde{F}(T) + O(\pi_L^{2n}) \,, 
  \]
  which is a contradiction. So when $\bar\alpha_1 = 1$, $\phi$ is the
  identity on~$D_P$.
  
  Otherwise, $\alpha_1 - 1$ is a unit, and then $F_\phi(T) = T$ has
  a unique solution in $\Fp_L$, by a standard Newton polygon
  argument. It follows that $\phi$ has a unique fixed 
  point on~$D_P(L)$. Since this holds for all fields~$L$ as above, 
  there is also a unique fixed point on~$D_P$.
\end{proof}

This result, a special case of which is the similar statement
on torsion points on elliptic curves
(see~\cite{SilvermanBook}, Prop.~VII.4.1), has the following
consequence.

\begin{Proposition} \label{PropImageOfI}
  Let $\alpha \in H^1(K, \Gamma)$ be a cohomology class, represented
  by a cocycle $\xi$. Then the restriction $\xi|_{I_K}$
  is a homomorphism $I_K \To \Gamma$, and its image is a cyclic
  subgroup of~$\Gamma$, whose conjugacy class only depends on~$\alpha$.
\end{Proposition}
\begin{proof}
  Since by Prop.~\ref{PropIActsTriv} $I_K$ acts trivially on~$\Gamma$,
  the cocycle condition for~$\xi$ on~$I_K$ simply means that 
  $\xi|_{I_K}$ is a homomorphism. Since $p \nmid \#\Gamma$, this
  homomorphism has to factor through a (finite) quotient of~$I_K$
  of order prime to~$p$. All such quotients are cyclic. The last
  statement follows from the definition of cohomology classes
  in~$H^1(I_K, \Gamma)$ (and again the fact that $I_K$ acts trivially).
\end{proof}

For $\xi \in Z^1(K, \Gamma)$, we define its {\em type} $t(\xi)$ to be
the cyclic subgroup
\[ t(\xi) = \xi(I_K) \subset \Gamma \,. \]

With these notations, we can state the main theorem in its local
version.

\begin{Theorem} \label{ThmMainLocal}
  Let $C/K$ be a smooth projective curve of genus $g \ge 1$, such
  that $C$ has good reduction. Choose $\Gamma$ as above, and
  pick a divisor class~$\CD$ of positive degree that is fixed
  under~$\Gamma$ in order to define the map $\phi : C \To J$.
  
  Let $\xi \in Z^1(K, \Gamma)$ be a cocycle whose cohomology class
  is ramified, and let $C_\xi$ be the corresponding twist of~$C$,
  with map $\phi_\xi : C_\xi \To J_\xi$. Assume that 
  \[ p > \#t(\xi)\,v(p) + f_{C/k}(r) + 1 \]
  for some $0 \le r < g$. Let $K'/K$ be an unramified extension, and let
  \[ F_\xi = \{P \in C_\xi(\bar{K}) \mid P \text{\ fixed under $t(\xi)$}\} \,.
  \]
  Pick a subgroup~$G$ of $J_\xi(K')$ of rank~$r$ and set
  \[ T = \{P \in C_\xi(K') \mid \phi_\xi(P) \in G\} \,. \]
  If $F_\xi$ is empty, then $T$ is empty, and otherwise
  \[ \#T \le f_{C/k}(r) + \#(T \cap F_\xi) 
              + \#(F_\xi \setminus (T \cup C_{\xi,\tors})) \,.
  \]
\end{Theorem}
                      
The proof will be given in section~\ref{SectProof} below.


\section{The main result --- global version}

We are now back in the global situation, where $C/\CK$ is a smooth
projective curve over a number field.

We let $\Gamma \subset \Aut_{\bar{\CK}}(C)$ be a finite $\CK$-defined
subgroup of the automorphism group of~$C$ and let 
\[ C^\triv
   = \{ P \in C \mid \gamma(P) = P
                     \text{\ for some\ } 1 \neq \gamma \in \Gamma \}
\]
be the set of {\em $\Gamma$-trivial points} on~$C$; $C^\triv$ is a finite (possibly empty) subset of~$C$. For the ``basepoint'' of the map of~$C$
into its Jacobian, we will always choose a divisor class of positive
degree that is invariant under~$\Gamma$.

We let $m(\Gamma)$ denote the maximal order of an element of~$\Gamma$
that fixes at least one point on~$C$.

Recall the notation $\tilde{f}_C$ and~$S_C$ from section~\ref{SectGeom}.

\begin{Theorem} \label{ThmTwistMaster}
  Let $C/\CK$ be a smooth projective curve of genus~$g$, and let $\Gamma$
  be a $\CK$-defined finite subgroup of the automorphism group of~$C$.
  Let $C_\xi$
  be a $\Gamma$-twist of~$C$ such that $\xi$ is ramified at some place~$v$
  of~$\CK$ outside~$S_C$.
  Assume that the residue characteristic~$p$ of~$v$ satisfies
  \[ p > m(\Gamma)\,e_v + \tilde{f}_C(r) + 1 \]
  for some $r < g$, where $e_v$ is the ramification index of~$\CK_v$ 
  over~$\Q_p$.
  Let $G$ be a subgroup of~$J_\xi(\CK)$ of rank~$r$ and let
  $T = \{P \in C_\xi(\CK) \mid \phi_\xi(P) \in G\}$.
  Then 
  \[ \#T \le \tilde{f}_C(r) + \#(T \cap C_\xi^\triv)
                    + \#(C_\xi^\triv \setminus (T \cup C_{\xi,\tors}))
     \,. \]
\end{Theorem}
\begin{proof}
  This follows from Thm.~\ref{ThmMainLocal}, applied to $C$ over~$\CK_v$,
  since $F_\xi \subset C_\xi^\triv$ and $\#t(\xi) \le m(\Gamma)$.
\end{proof}

The bound can be read as ``$\tilde{f}_C(r)$, plus the trivial points in~$T$,
plus the trivial non-torsion points outside~$T$''. The last contribution
is somewhat annoying (the first can be considered as bounding the
non-trivial points in~$T$). Luckily, in many cases, all the trivial
points are torsion, and so this contribution vanishes. A sufficient
condition for this is that all the quotients $C/\langle \gamma \rangle$
have genus zero, for all $1 \neq \gamma \in \Gamma$ having at least
one fixed point on~$C$.

In cases where this last contribution cannot be shown to vanish,
it may be possible to improve the bound by noting that it can be
replaced by
\[ \max\{\#(C^{\langle \gamma \rangle} \setminus C_\tors)
          \mid 1 \neq \gamma \in \Gamma\} \,, 
\]
the maximal number of non-torsion fixed points of any nontrivial
automorphism in~$\Gamma$,
which bounds $\#(F_\xi \setminus (T \cup C_{\xi,\tors}))$ in
Thm.~\ref{ThmMainLocal}.

\medskip

Now let us proceed to prove the results given in the introduction.
To prove Thm.~\ref{IntroThmHC}, we set $\CK = \Q$ and $\Gamma = \mu_2$
acting on~$y$. Then $C^\triv$ consists of the Weierstrass points, which
map to zero in the Jacobian, i.e., $C^\triv \subset C_\tors$, and the
same holds for all the quadratic twists. Hence the last term in the
bound above vanishes. Now take $\Sigma \subset C_d(\Q)$ as in Thm.~\ref{IntroThmHC}
and let $G$ be the subgroup generated by~$\Sigma$ in the Jacobian.
The set $T$ in our theorem above then contains $\Sigma \cup \iota(\Sigma)$ and
the rational Weierstrass points of~$C_d$. We assume that $r < \#\Sigma$; then
\[ 2\#\Sigma + \#(T \cap C_d^\triv) \le \#T \le 2 r + \#(T \cap C_d^\triv) \,, \]
a contradiction. Therefore, we must have $r = \#\Sigma$. 
Note that the condition on ramification
reduces to $p\,|\,d$ with $p > 2\#\Sigma + 1$ and such that $C$ has good
reduction at~$p$\,; the latter condition is implied by $p$ not dividing
the discriminant of~$f$.

The proof of Thm.~\ref{IntroThmHCVariant} is very similar. Here we
take for $G$ the Mordell-Weil group $J_d(\Q)$; then $T = C_d(\Q)$
in Thm.~\ref{ThmTwistMaster}, which then says that
\[ \#C_d(\Q) \le 2r + \#\{\text{Weierstrass points in~$C_d(\Q)$}\} \,. \]

To prove Thm.~\ref{IntroThmThue}, we consider the smooth plane projective
curve 
\[ C : Z^n = F(X, Y) \,. \] 
We take $\Gamma = \mu_n$, acting on~$Z$,
and we let $G$ be the Mordell-Weil group in the Jacobian of the twist
$C_h : h\,Z^n = F(X, Y)$. If its rank $r$ is at most $n-3$, then $f_C(r) = r$,
and we get the same value for the reduction mod~$p$ when $p$ does not
divide the discriminant of~$F$.
Furthermore, the trivial points all have $Z = 0$ (and they all belong
to~$C_\tors$), so we do not see them in the affine equation~\eqref{EqnThue}.
Hence only the first term in the bound remains, which is~$r$. The ramification
conditions again boil down to $p\,|\,h$ with $p > n + r + 1$ and
$p \nmid \disc(F)$. The more general bound for larger~$r$ follows
similarly, using that $f_C(r) \le 2r$. Note that $g = \half(n-1)(n-2)$.

The next result, Thm.~\ref{IntroThmPointsOnTwists}, follows directly
from Thm.~\ref{ThmTwistMaster} above, noting that the ramification condition is
violated only for a finite set of cocycle classes. Here we pick a set~$S_C$
that is large enough so that we can choose $\tilde{f}_C = f_C$.

Another easy application is the following.

\begin{Proposition}
  Let $\ell \ge 9$ be an odd integer. Then there are infinitely many
  $2\ell$th power free integers~$A$ such that $r(C_{\ell,A}) \ge 4$,
  where $C_{\ell, A} : y^2 = x^\ell + A$.
\end{Proposition}
\begin{proof}
  We take $C = C_{\ell,1}$ and $\Gamma = \mu_2 \times \mu_\ell$. Then we
  can deduce that for all but finitely many~$A$, a set of $n$ points
  in~$C_{\ell,A}(\Q)$ with nonvanishing $x$-coordinate and
  positive $y$-coordinate
  will generate a subgroup of rank~$n$ in the Mordell-Weil group,
  provided $n \le g = (\ell-1)/2 \ge 4$. Since the family
  \[ Y^2 + Y = X^\ell + t^{2\ell} \]
  provides infinitely many curves $C_{\ell,A}$ that have at least four such 
  points (given by $X \in \{-t, t, -t^2, t^4\}$), the claim follows.
\end{proof}

The statement is also true for $\ell = 5$ and~$7$, as can be shown by
more direct methods.


\section{The Chabauty-Coleman machine}

The Chabauty-Coleman method consists in bounding the number of 
$\CK$-rational points of~$C$, where $\CK$ is a number field,
in terms of the number of zeros of certain differentials of $C/K$, where
$K = \CK_v$ for some finite place $v$ of~$\CK$. We will extend the method
and allow $K$ to be some finite extension of~$\CK_v$.

The basic result of this method is a bound for the number of $\CK_v$-points
mapping into a subgroup of rank $r < g$ in the Jacobian, which takes
the form (see Thm.~\ref{ThmBound} below)
\[ \#\CX + f_{C/k_v}(r) + \Delta_v(\#\CX, f_{C/k_v}(r)) \,. \]
The first term is the number of residue classes the points are sitting in,
the second term comes from the zeros of differentials and was up to now
taken to be $2g-2$ in published applications of the method. We improve this
part of the bound by replacing it with~$f_{C/k_v}(r)$. The third term is a
contribution that has to be put in to account for the possible presence
of small denominators divisible by~$p$ in the logarithm series. It vanishes
if $p$ is large enough. On the other hand, the first term usually grows
with~$p$. The new trick we use to obtain our main result lets us keep
$\#\CX$ small while $p$ is large.

Since we need the precise statement with the improvement we introduce,
we provide proofs of the relevant facts.

Let us first set some notation. As before, $K$ will be a $p$-adic local field 
with normalised valuation~$v$, uniformiser~$\pi$
and residue class field~$k$ (of characteristic~$p$). We will denote
the ring of integers of~$K$ by~$\CO$.

Sometimes we will consider a finite field extension $L/K$; then we
denote the objects associated to~$L$ by $v_L$, $k_L$, $\CO_L$ etc.

We let $e = v(p)$ be the ramification index of $K/\Q_p$ and define
\begin{align*}
   \delta(v, n) 
     &= \max\{ d \ge 0 \mid n+d+1 - v(n+d+1) \le n+1 - v(n+1) \} \\
     &= \max\{ d \ge 0 \mid e\,v_p(n+1) + d \le e\,v_p(n+d+1) \} \,.
\end{align*}
For $s, r \ge 0$, let
\[ \Delta_v(s, r)
     = \max\Bigl\{ \sum_{j=1}^s \delta(v,m_j)
                    \Bigm| \sum_{j=1}^s m_j \le r \Bigr\} \,.
\]
Note that $\Delta_v$ is obviously an increasing function in both arguments.

We need to bound $\delta$ and~$\Delta$ from above.

\begin{Lemma} \label{Delta1}
  If $p > e + 1$, then $\delta(v, n) \le e\,\lfloor n/(p-e-1) \rfloor$.
  In particular, if $p > n + e + 1$, then $\delta(v, n) = 0$.
\end{Lemma}
\begin{proof}
  If $\delta(v, n) = d$, then $e v_p(n+d+1) \ge e v_p(n+1) + d \ge d$.
  This implies that $p^{\ceil{d/e}}$ divides 
  $n+d+1$, so $p^{\ceil{d/e}} \le n + d + 1$. Hence we always have
  \[ \delta(v, n) \le e\,\max\{ d \mid p^d \le n + ed + 1 \} \,. \]
  Now suppose $p \ge e + 2$. Then by an easy induction, we see that
  $p^d - ed - 1 \ge (p-e-1) d$ for all $d \ge 0$, hence
  $d \le n/(p-e-1)$ for all~$d$ in the set above.
\end{proof}

\begin{Lemma} \label{Delta2} 
  If $p > e + 1$, we have
  \[ \Delta_v(s, r) \le e\,\lfloor r/(p-e-1) \rfloor \,. \]
  In particular, if $p > r + e + 1$, then $\Delta_v(s, r) = 0$.
\end{Lemma}
\begin{proof}
  This follows immediately from Lemma~\ref{Delta1} and the definitions.
\end{proof}

In the following, we will assume that $C/K$ has good reduction. 
We will denote by $\CC$ a smooth model of $C$ over~$\CO$. 
Let $\CC_\s$ denote the special fiber of~$\CC$.
The preimage of a point $P \in \CC_\s(\bar{k})$ under
the reduction map $\rho : C(\bar{K}) \To \CC_\s(\bar{k})$ 
is called a {\em residue class} and denoted by~$D_P$. 
We write $D_P(K)$ for $D_P \cap C(K)$.

We write $\Omega(C/K)$ for the 
$K$-vector space of global regular ($K$-rational) differentials on~$C$.
This space has dimension~$g$.

There is a pairing
\[ \Omega(C/K) \times J(K) \To K \]
with trivial left kernel and right kernel $J(K)_{\tors}$; it is
given by a logarithm map, see \cite{McCallum94} or~\cite{Wetherell}. 
So when the rank $r(C)$
of~$J(\CK)$ is less than $g$, the dimension of~$\Omega(C/K)$,
then there is a subspace $\Lambda$ of~$\Omega(C/K)$ of codimension 
at most~$r(C)$ that annihilates~$J(\CK)$.

We will use a more general setup here. Let $G \subset J(K)$ be a
subgroup of torsion-free rank $r < g$ (i.e., $r = \dim_\Q G \tensor_\Z \Q$).
Furthermore, let $\bar{G}$ be the saturation of~$G$, i.e.,
\[ \bar{G} = \{ P \in J(K) \mid nP \in G \text{\ for some\ } n > 0 \} \,. \]
We define
\[ \Lambda = \Lambda(G)
           = \{ \omega \in \Omega(C/K) \mid \omega \text{\ kills\ } G \} \,.
\]
Then we obviously have $\codim \Lambda \le r$, i.e., $\dim \Lambda \ge g-r$.
Note that $\Lambda$ kills not only~$G$, but also~$\bar{G}$.

Let $X$ be the set of points in~$C(K)$ that map into $\bar{G}$
under the fixed map of~$C$ into~$J$. For example, we could have $G = J(\CK)$,
then $X \supset C(\CK) \cup C(K)_\tors$. We let $\CX$ denote the image of~$X$
in~$\CC_\s(k)$.

Let $0 \neq \omega \in \Omega(C/K)$. Then (since $\Omega(\CC/\CO)$ is
a lattice in~$\Omega(C/K)$) there is a multiple of~$\omega$
reducing to a nonzero differential $\bar\omega \in \Omega(\CC_\s/k)$.
If $P \in \CC_\s(k)$ is a point, 
we denote by $n(\omega,P) = v_P(\bar\omega)$ the order of vanishing 
of~$\bar\omega$ in~$P$.
We write $\nu(P) = \#(D_P \cap X)$. 

Then we have the following result.

\begin{Proposition} \label{PropRC}
  Let $0 \neq \omega \in \Lambda$, and let $P \in \CC_s(k)$. Then
  \[ \nu(P) \le 1 + n(\omega, P) + \delta(v, n(\omega, P)) \,. \]
\end{Proposition}
\begin{proof}
  Without loss of generality, $\omega$ itself reduces to~$\bar\omega$.
  For simplicity, write $n$ for $n(\omega, P)$.
  Choose a uniformizer $t$ at a point in~$D_P(K)$; since 
  \[ \bar\omega = (u\,\bar{t}^n + \text{higher order terms})\,d\bar{t} \]
  with $u \in k^\times$,
  $\omega$ has an expansion with coefficients in~$\CO$,
  \[ \omega = (a_0 + a_1\,t + a_2\,t^2 + \dots)\,dt \,, \]
  where $a_0$, $a_1$, \dots, $a_{n-1}$ have positive valuation and
  $a_n$ is a $v$-adic unit. The logarithm corresponding to~$\omega$ is
  then given on~$D_P(K)$ by
  \[ \lambda_\omega(Q)
       = c + a_0\,\pi\,T + \dots + \frac{a_m}{m+1}\,\pi^{m+1}\,T^{m+1} + \dots
       = \lambda(T)
       \,,
  \]
  where $t(Q) = \pi\,T$, $T \in \CO$, and $c$ is a constant of
  integration. By a standard Newton polygon argument, $\lambda_\omega$ has
  at most $n + 1 + \delta(v, n)$ zeros on~$D_P(K)$, since these zeros
  correspond to integral zeros of the power series~$\lambda$.

  Since $X \cap D_P$ is contained in this set of zeros (for 
  $\lambda_\omega$ kills~$X$ by definition), the claim follows.
\end{proof}

This result prompts the following definitions.
\begin{align*}
  n(\Lambda, P) &= \min \{ n(\omega, P) \mid 0 \neq \omega \in \Lambda \} 
    \text{\qquad and} \\
  N(\Lambda, C/K) &= \sum_P n(\Lambda, P) \,,
\end{align*}
where the sum extends over the points $P \in \CC_\s(k)$.

In their 
paper~\cite{LorenziniTucker}, p.~59, Lorenzini and Tucker ask whether it
is possible to sharpen the trivial bound $N(\Lambda, C/K) \le 2g-2$ to get
$N(\Lambda, C/K) \le 2 \codim \Lambda$. 
We can give an affirmative (and even better) answer.
Let $f_{C/k}$ be defined as in section~\ref{SectGeom}.

\begin{Theorem} \label{ThmNLambda}
  Let $C/K$ be a smooth projective curve of genus~$g$, 
  and let $0 \neq \Lambda$ be a $K$-linear 
  subspace of~$\Omega(C/K)$. If $C$ has good reduction, then
  \[ N(\Lambda, C/K) \le f_{C/k}(\codim \Lambda) \le 2 \codim \Lambda \,. \]
\end{Theorem}
\begin{proof}
  Because of good reduction, there is a well-defined reduction map
  \[ \rho : \BP(\Omega(C/K)) \To \BP(\Omega(\CC_s/k)) \,, \]
  which preserves
  dimensions of subspaces. Let $\bar\Lambda$ be the linear subspace
  corresponding to the image of~$\BP(\Lambda)$; it has
  dimension $\dim \Lambda$. For any $\omega \in \BP(\Omega(C/K))$,
  we have $n(\omega, P) = v_P(\rho(\omega))$. Let $D$ be the effective divisor
  \[ D = \sum_P n(\Lambda, P) \cdot P \]
  on $\CC_\s$; then $N(\Lambda, C/K) = \deg D \le 2g-2$ and
  $\omega \in \bar\Lambda$ implies $(\omega) \ge D$.
  So 
  \[ N(\Lambda, C/K)
       \le \max\{\deg D \mid D \ge 0, \dim \Omega(D) \ge \dim \Lambda\}
       = f_{C/k}(\codim \Lambda)\,.
  \]
\end{proof}

\begin{Remark}
  This result is still true when $C$ is hyperelliptic and of {\em bad}
  reduction. However, to formulate the method in the case of bad reduction
  requires the use of a minimal proper regular model, and since we do not
  need this case here, we refrain from giving the details.
\end{Remark}

\medskip

We now have the following ``master theorem'' for the
Chabauty-Coleman method.

\begin{Theorem} \label{ThmBound}
  Let $G \subset J(K)$ be a subgroup of rank $r < g$, let
  \[ X = \{P \in C(K) \mid \phi(P) \in \bar{G}\} \,, \] 
  and let $\CX$ be the image of~$X$ in~$\CC_\s$.
  Then we have the bound
  \[ \#X \le \#\CX + f_{C/k}(r) + \Delta_v(\#\CX, f_{C/k}(r))
     \,.
  \]
  Furthermore, when $p > f_{C/k}(r) + e + 1$, we have
  \[ \#X \le \#\CX + f_{C/k}(r) \,. \]
\end{Theorem}
\begin{proof}
  Sum the bounds of Prop.~\ref{PropRC} over the residue classes corresponding
  to points in~$\CX$ and use Thm.~\ref{ThmNLambda} and Lemma~\ref{Delta2}.
\end{proof}

As a corollary, we get a refinement of Coleman's bound~\cite{Coleman}.

\begin{Corollary}
  Let $\CK$ be a number field, $C/\CK$ a curve, and suppose that
  $C$ has good reduction at the finite place $v$ of~$\CK$. 
  Let $K = \CK_v$, and let $\CC/\CO$ be a model of~$C$ with good reduction. Then
  if the rank $r(C)$ of~$J(\CK)$ is less than the genus of~$C$, we have
  \begin{align*}
    \#C(\CK) &\le \#\CC(k) + f_{C/k_v}(r(C)) + \Delta_v(\#\CC(k), f_C(r(C)))\\
             &\le \#\CC(k) + 2\,r(C) + \Delta_v(\#\CC(k), 2\,r(C)) \,.
  \end{align*}
  If in addition $p > f_{C/k_v}(r(C)) + e_v + 1$, where $p$ is the residue
  characteristic of~$v$, then we have
  \[ \#C(\CK) \le \#\CC(k) + f_{C/k_v}(r(C)) \le \#\CC(k) + 2\,r(C) \,. \]
\end{Corollary}
\begin{proof}
  We choose $G = J(\CK) \subset J(K)$. The set~$X$ then contains~$C(\CK)$, and
  using Thm.~\ref{ThmBound}, we get the bound as given, noting that
  $f_{C/k_v}(r(C)) \le 2 r(C)$ by Lemma~\ref{LemmaFBound}.
\end{proof}


\section{Proof of the main result} \label{SectProof}

Now we want to apply this machine to prove our main result.
We continue to use the notations set in the previous sections.
In particular, we will always assume that $v(\#\Gamma) = 0$, i.e.,
that $p$ does not divide the order of~$\Gamma$, where $p$ is the
residue characteristic of~$K$.

Recall that a cocycle class in $H^1(K, \Gamma)$ is called
{\em ramified} when it has nontrivial image in $H^1(K^\unr, \Gamma)$,
where $K^\unr$ is the maximal unramified extension of~$K$.

\begin{Proposition} \label{PropTrivialImage}
  Let $C/K$ be a curve with good reduction, and
  let $C_\xi$ be a $\Gamma$-twist of~$C$ such that $\xi$ is ramified. 
  Let $L/K$ be a finite extension such that $C$ and $C_\xi$ become
  isomorphic over~$L$, and let $\CC$ be a model of~$C$ over~$\CO$
  that has good reduction. Then the image of~$C_\xi(K)$ in~$\CC_\s(k_L)$
  consists of points fixed under~$t(\xi)$, the type of~$\xi$.
\end{Proposition}
\begin{proof}
  Let $\varphi : C_\xi/L \To C/L$ be an isomorphism, and let $(\xi_\sigma)$
  be the associated cocycle $\xi_\sigma = \varphi^\sigma \varphi^{-1}$
  taking values in~$\Gamma$.
  We have $t(\xi) = \xi(I_K) \neq \{1\}$, and for all 
  $\gamma = \xi_\sigma \in t(\xi)$ (for some $\sigma \in I_K$) and all points
  $P \in C_\xi(K)$, we have (indicating images in~$\CC_\s$ by
  putting a bar above the point)
  \[ \gamma(\overline{\varphi(P)})
       = \overline{\gamma(\varphi(P))}
       = \overline{\varphi^\sigma(P)}
       = \overline{\varphi(P^{\sigma^{-1}})^\sigma}
       = \overline{\varphi(P^{\sigma^{-1}})}
       = \overline{\varphi(P)} \,,
  \]
  so $\overline{\varphi(P)}$ is fixed by~$\gamma$. 
  The last two equalities use that
  $\sigma$ acts trivially both on~$\CC_\s$ and on~$P$.
\end{proof}

\begin{Corollary} \label{CorEmpty}
  Let $C/K$ be a curve with good reduction, and
  let $C_\xi$ be a $\Gamma$-twist of~$C$ such that $\xi$ is ramified. 
  Assume that $v(\#\Gamma) = 0$ and that $C^\triv(\bar{K}) = \emptyset$ (or just
  $C^{t(\xi)}(\bar{K}) = \emptyset$). Then $C_\xi(K) = \emptyset$.
\end{Corollary}
\begin{proof}
  Let $L$ be as in Prop.~\ref{PropTrivialImage}.
  The condition $v(\#\Gamma) = 0$, together with good reduction, implies
  that $C^{t(\xi)}(L) \to \CC_\s^{t(\xi)}(k_L)$ is a bijection. Since by assumption,
  $C^{t(\xi)}(L)$ is empty, so is $\CC_\s^{t(\xi)}(k_L)$. Since by 
  Prop.~\ref{PropTrivialImage}, $C_\xi(K)$ maps into~$\CC_\s^{t(\xi)}(k_L)$,
  $C_\xi(K)$ must be empty.
\end{proof}

If we apply this to curves over number fields, we get the following
well-known result. (This goes back to Chevalley and Weil~\cite{ChevalleyWeil}.)

\begin{Theorem}
  Let $D \To C$ be an unramified cover of curves over a number
  field~$\CK$ that is geometrically Galois with Galois group~$\Gamma$. 
  Then for all but finitely many
  twists $D' \To C$ of this cover, $D'(\CK)$ is empty. More strongly,
  only finitely many twists have points everywhere locally.
\end{Theorem}
\begin{proof}
  All but finitely many twists are ramified at some place $v$ with
  $v(\#\Gamma) = 0$ and such that $D$ has good reduction at~$v$. 
  For such a twist, Cor.~\ref{CorEmpty} shows that
  already $D'(\CK_v)$ is empty. Note that the twists of the cover are
  exactly the $\Gamma$-twists of~$D$ and that $D^\triv$ is empty, since
  we assume the cover to be unramified.
\end{proof}

\medskip

Now let us proceed to prove Thm·~\ref{ThmMainLocal}.
Note that the assumptions are preserved if we replace $K$ by an unramified
extension. We can therefore assume that $K' = K$.
There is a finite extension $L/K$ such that
$C$ and $C_\xi$ are isomorphic over~$L$, which
can be taken to be the fixed field of 
$\{\sigma \in G_K \mid \xi_\sigma = 1\}$. Then $e_{L/K} = \#t(\xi)$,
and so $v_L(p) = \#t(\xi) v(p)$.
We apply Thm.~\ref{ThmBound}
with the field~$L$ and the group $G \subset J_\xi(K) \inj J(L)$.
As before, let $\CC$ be a model of~$C$ over~$\CO$ with good reduction.
Then $X$ is the set of points in~$C(L)$ mapping into the saturation~$\bar{G}$
of~$G$, and $\CX$ is the image of~$X$ in~$\CC_\s(k_L)$.
By Prop.~\ref{PropTrivialImage}, the set $\CX$ consists of fixed
points of~$t(\xi)$ in~$\CC_\s(k_L)$. Now Thm.~\ref{ThmBound} says
\[ \#X \le \#\CX + f_{C/k}(r) + \Delta_w(\#\CX, f_{C/k}(r)) \,, \]
where $w = v_L$ is the normalised valuation of~$L$. 
Now $e_L = \#t(\xi) e_K = \#t(\xi) v(p)$, so we 
have $p > f_{C/k}(r) + e_L + 1$, which gives the better bound
\[ \#X \le \#\CX + f_{C/k}(r) \,. \]
In the following, we will identify $T$ with its image in~$C(L)$, so
$T \subset X$. We have $v(\#\Gamma) = 0$. This implies
that each point $P \in \CC_\s(k_L)^{t(\xi)}$ lifts to a (unique) point
$\tilde{P} \in D_P \cap C^{t(\xi)}(L)$. Such a point will belong to~$X$
if it belongs to~$T$ or if it is torsion. Hence
\[ \#\CX \le \#\CC_\s^{t(\xi)}(k_L) = \#C^{t(\xi)}(L)
         \le \#(C^{t(\xi)}(L) \setminus (T \cup C_\tors)) 
              + \#(X \cap C^{t(\xi)}(L))
   \,. 
\]
Therefore
\begin{align*}
  \#T &\le \#(X \setminus C^{t(\xi)}(L)) + \#(T \cap C^{t(\xi)}(L)) \\
      &= \#X - \#(X \cap C^{t(\xi)}(L)) + \#(T \cap C^{t(\xi)}(L)) \\
      &\le f_{C/k}(r) + \#\CX - \#(X \cap C^{t(\xi)}(L)) 
          + \#(T \cap C^{t(\xi)}(L)) \\
      &\le f_{C/k}(r) + \#(T \cap C^{t(\xi)}(L))
          + \#(C^{t(\xi)}(L) \setminus (T \cup C_\tors)) \\
      &\le f_{C/k}(r) + \#(T \cap F_\xi) 
          + \#(F_\xi \setminus (T \cup C_{\xi,\tors}))
\end{align*}
as it was stated in Thm.~\ref{ThmMainLocal}. Note that in the last line,
we have used the identification of $C_\xi(L)$ and~$C(L)$ given by the
isomorphism of $C_\xi$ and~$C$ over~$L$ in order to transfer the result
back to~$C_\xi$.



\begin{thebibliography}{CHMMM}

\frenchspacing
\renewcommand{\baselinestretch}{1}\small

\bibitem[CHM97]{CHM}
  {\sc L. Caporaso, J. Harris} and {\sc B. Mazur:}
  {\it Uniformity of rational points}, J. Am. Math. Soc. {\bf 10}, 1--35
  (1997).
\bibitem[CC55]{CartanChevalley}
  {\sc H. Cartan} and {\sc C. Chevalley:}
  {\it G\'eom\'etrie alg\'ebrique},
  S\'eminaire Cartan-Chevalley, Secr\'etariat Math., Paris (1955/56).
\bibitem[CW32]{ChevalleyWeil}
  {\sc C. Chevalley, A. Weil:} {\it Un th\'eor\`eme d'arithm\'etique
   sur les courbes alg\'ebriques,} Comptes Rendus Hebdomadaires des
   S\'eances de l'Acad. des Sci., Paris {\bf 195}, 570--572 (1932).
\bibitem[Col85]{Coleman}
  {\sc R.F. Coleman:} {\it Effective Chabauty}, Duke Math. J. {\bf 52},
  765--770 (1985).
\bibitem[Elk99]{Elkies}
  {\sc N.D. Elkies:}
  {\it The Klein Quartic in number theory}, in:
  {\it The eightfold way}, Cambridge University Press, 1999, pp.~51--101.
\bibitem[Har77]{Hartshorne}
  {\sc R. Hartshorne:} {\it Algebraic Geometry},
  Springer-Verlag, New York etc. (1977).
\bibitem[LT02]{LorenziniTucker}
  {\sc D. Lorenzini} and {\sc T.J. Tucker:} {\it Thue equations and the
  method of Chabauty-Coleman}, Invent. Math {\bf 148}, 47--77 (2002).
\bibitem[McC94]{McCallum94}
  {\sc W.G. McCallum:} {\it On the method of Coleman and Chabauty},
  Math. Ann. {\bf 299}, 565--596 (1994).
\bibitem[Pac97]{Pacelli}
  {\sc P.L. Pacelli:}
  {\it Uniform boundedness for rational points},
  Duke Math. J. {\bf 88}, 77--102 (1997).
\bibitem[Ray83]{Raynaud}
  {\sc M. Raynaud:}
  {\it Courbes sur une vari\'et\'e ab\'elienne et points de torsion},
  Invent. Math. {\bf 71}, 207--233 (1983).
\bibitem[Ser97]{Serre}
  {\sc J.-P. Serre:}
  {\it Galois cohomology},
  Springer-Verlag, Berlin-Heidelberg-New York (1997).
\bibitem[Sil86]{SilvermanBook}
  {\sc J.H.\ Silverman:} {\it The arithmetic of elliptic curves},
  Grad.\ Texts in Math.\ 106, Springer-Verlag, New York (1986).
\bibitem[Sil93]{Silverman93}
  {\sc J.H. Silverman:} {\it A uniform bound for rational points on twists
  of a given curve}, J. London Math. Soc. {\bf 47}, 385--394 (1993).
\bibitem[Sil94]{SilvermanBook2}
  {\sc J.H.\ Silverman:} 
  {\it Advanced topics in the arithmetic of elliptic curves},
  Grad.\ Texts in Math.\ 151, Springer-Verlag, New York (1994).
\bibitem[Sto06]{Stoll06}
  {\sc M. Stoll:} 
  {\it On the number of rational squares at fixed distance from a fifth power,}
  Preprint (2006).
\bibitem[Wet97]{Wetherell}
  {\sc J.L. Wetherell:} {\it Bounding the number of rational points on
  certain curves of high rank}, Ph.D. thesis, University of California
  (1997).

\end{thebibliography}
\end{document}